\documentclass{gtart_h}  

\def\ifplaintex{\expandafter\ifx\csname documentclass\endcsname\relax}

\def\ifplaintex{\expandafter\ifx\csname documentclass\endcsname\relax}

\ifplaintex 
\else
\fi

\expandafter\ifx\csname epsfbox\endcsname\relax\input epsf\fi

\def\gt{{\mathsurround=0pt\it $\cal G\mskip-2mu$eometry \&\ 
$\cal T\!\!$opology}}        

\def\gtp{{\mathsurround=0pt\it $\cal G\mskip-2mu$eometry \&\ 
$\cal T\!\!$opology $\cal P\!$ublications}}  

\def\lognumber#1{\def\thelognumber{#1}}
\def\volumenumber#1{\def\thevolumenumber{#1}}
\def\papernumber#1{\def\thepapernumber{#1}}
\def\volumeyear#1{\def\thevolumeyear{#1}}

\def\pagenumbers#1#2{\def\startpage{#1}\def\finishpage{#2}}
\def\published#1{\def\publishdate{#1}}
\def\proposed#1{\def\theproposer{#1}}
\def\seconded#1{\def\theseconders{#1}}
\def\received#1{\def\receiveddate{#1}}

\def\accepted#1{\def\accepteddate{#1}}
\def\asciititle#1{\def\theasciititle{#1}}

\def\asciiaddress#1{\def\theasciiaddress{#1}}

\long\def\asciiabstract#1{\long\def\theasciiabstract{#1}}
\def\asciikeywords#1{\def\theasciikeywords{#1}}

\let\\\par\let\thelognumber\relax
\let\thevolumenumber\relax\let\thepapernumber\relax
\let\thevolumeyear\relax\let\thesamplenumber\relax\let\startpage\relax
\let\finishpage\relax\let\publishdate\relax\let\receiveddate\relax
\let\reviseddate\relax\let\accepteddate\relax\let\theasciititle\relax
\let\theasciiauthors\relax\let\theasciiaddress\relax
\let\theasciiabstract\relax\let\theasciikeywords\relax
\let\theasciiemail\relax\let\theshortauthors\relax\let\theshorttitle\relax

\long\def\maketitlep{   

\count0=\startpage

\gt\hfill      
\hbox to 77pt{\vbox to 0pt{\vglue -15pt\epsfbox{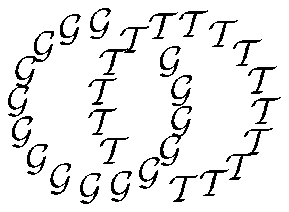}\vss}\hss}
\break
{\small\ifx\thesamplenumber\relax 
Volume \else Sample
\fi\thevolumenumber\ (\thevolumeyear)
\startpage--\finishpage\nl
Published: \publishdate}
\vglue 0.5truein plus 0.4fil minus 0.1truein

{\parskip=0pt\leftskip 0pt plus 1fil\def\\{\par\smallskip}{\ifplaintex\large
\else\Large\fi\bf\thetitle}\par\medskip}   

\vglue 0pt plus 0.1fil 

{\parskip=0pt\leftskip 0pt plus 1fil\def\\{\par}{\sc\theauthors}
\par\medskip}

\vglue 0pt plus 0.1fil 

{\small\parskip=0pt\let\newline\\
{\leftskip 0pt plus 1fil\def\\{\par}{\sl\theaddress}\par}
\expandafter\ifx\theemail\relax    
\relax\else\vglue 5pt plus 0.02fil minus 2pt\def\\{\stdspace{\rm 
and}\stdspace} 
\cl{Email:\stdspace\tt\theemail}\fi
\ifx\theurl\relax                  
\relax\else\vglue 5pt plus 0.02fil minus 2pt\def\\{\stdspace{\rm 
and}\stdspace}
\cl{URL:\stdspace\tt\theurl}\fi\par}

\vglue 7pt plus 0.3fil minus 3pt

{\bf Abstract}
\vglue 5pt plus 0.1fil minus 2pt

\theabstract

\vglue 7pt plus 0.3fil minus 3pt

{\bf AMS Classification numbers}\quad Primary:\quad \theprimaryclass

Secondary:\quad \thesecondaryclass

\vglue 5pt plus 0.3fil minus 2pt

{\bf Keywords:}\quad \thekeywords

\vglue 10pt plus 0.5fil minus 5pt

{\small  Proposed: \theproposer\hfill Received: \receiveddate\nl
Seconded: \theseconders\hfill 
\ifx\reviseddate\relax                         
Accepted: \accepteddate                        
\else
Revised: \reviseddate                          
\fi}
\eject
}       

\font\phead=cmsl9 scaled 950
\font=cmsl9 scaled 1050
\font\pnum=cmbx10 scaled 913
\font\lnum=cmbx10 
\font\pfoot=cmsl9 scaled 950
\font=cmsl9 scaled 1050
\ifplaintex
\headline{\vbox to 0pt{\vskip -4.5mm\line{\small\phead\ifnum
\count0=\startpage ISSN 1364-0380 (on line)
1465-3060 (printed) \hfill {\pnum\folio}\else\ifodd\count0\def\\{ }%
\ifx\theshorttitle\relax\thetitle\else\theshorttitle\fi\hfill{\pnum\folio}
\else\def\\{ and }{\pnum\folio}\hfill\ifx\theshortauthors\relax\theauthors
\else\theshortauthors\fi\fi\fi}\vss}}
\footline{\vbox to 0pt{\vglue 0mm\line{\small\pfoot\ifnum\count0=\startpage
\copyright\ \gtp\hfill\else
\gt, Volume \thevolumenumber\ (\thevolumeyear)\hfill\fi}\vss
}}
\else
\makeatletter
\def\@oddhead{{\small\lhead\ifnum\count0=\startpage ISSN 1364-0380 (on line)
1465-3060 (printed) \hfill {\lnum\number\count0}\else\ifodd\count0
\def\\{ }\ifx\theshorttitle\relax \thetitle \else\theshorttitle\fi\hfill
{\lnum\number\count0}\else\def\\{ and }{\lnum\number\count0}
\hfill\ifx\theshortauthors\relax 
\theauthors\else\theshortauthors\fi\fi\fi}}\def\@evenhead{\@oddhead}
\def\@oddfoot{\small\lfoot\ifnum\count0=\startpage\copyright\ \gtp\hfill\else
\gt, Volume \thevolumenumber\ (\thevolumeyear)\hfill\fi}
\def\@evenfoot{\@oddfoot}
\makeatother
\fi

\newwrite\gtoutfile
\long\gdef\makeheadfile{  
{\def\\{, }\def\s{ }
\immediate\openout\gtoutfile head.xxx
\immediate\write\gtoutfile{Proxy-for: \ifx\theasciiauthors\relax
\theauthors\else\theasciiauthors\fi\s<\ifx\theasciiemail\relax\theemail\else\theasciiemail\fi>}
\immediate\write\gtoutfile{\noexpand\\}
\immediate\write\gtoutfile{Authors: \ifx\theasciiauthors\relax
\theauthors\else\theasciiauthors\fi}
{\def\\{ }\immediate\write\gtoutfile{Title: \ifx\theasciititle\relax
\thetitle\else\theasciititle\fi}}
\immediate\write\gtoutfile{Subj-class: GT or SG or MG etc}
\immediate\write\gtoutfile{MSC-class: \theprimaryclass\ifx\thesecondaryclass\relax\else, \thesecondaryclass\fi}
\immediate\write\gtoutfile{Journal-ref: Geom. Topol. \thevolumenumber
(\thevolumeyear) \startpage-\finishpage}
\immediate\write\gtoutfile{Comments: Published by Geometry and Topology at}
\immediate\write\gtoutfile{\s\s http://www.maths.warwick.ac.uk/gt/GTVol\thevolumenumber/paper\thepapernumber.abs.html}
\immediate\write\gtoutfile{\noexpand\\}
\immediate\write\gtoutfile{}
\ifx\theasciiabstract\relax
\immediate\write\gtoutfile{\theabstract}\else
\immediate\write\gtoutfile{\theasciiabstract}\fi
\immediate\write\gtoutfile{}
\immediate\write\gtoutfile{\noexpand\\}
\immediate\write\gtoutfile{}
\immediate\closeout\gtoutfile}}  

\def\maketitlepage{\maketitlep\makeheadfile}
\let\maketitle\maketitlepage

\lognumber{627}
\received{23 June 2005}
\volumenumber{9}\papernumber{38}\volumeyear{2005}
\pagenumbers{1677}{1687}   
\published{28 August 2005}
\accepted{4 August 2005}
\proposed{Peter Ozsv\'ath}
\seconded{Robion Kirby, Yasha Eliashberg}

\usepackage{amscd, amsmath, amssymb, graphicx, psfrag}
\usepackage[all]{xy}

\def\psfraga <#1,#2> #3#4{%
\psfrag {#3}{\smash{\rlap{\kern #1 \raise #2\hbox{#4}}}}}

\def\figref#1{\hyperlink{#1anchor}{Figure~\ref*{#1}}}
\def\anchor#1{\noindent\hypertarget{#1anchor}{\smash{$\phantom{99}$}}\newline}

\newcommand{\R}{\mathbb R}
\newcommand{\N}{\mathbb N}
\newcommand{\Q}{\mathbb Q}
\newcommand{\Z}{\mathbb Z} 

\theoremstyle{plain}
\newtheorem{theorem}{Theorem}[section]
\newtheorem{prop}[theorem]{Proposition}
\newtheorem{lemma}[theorem]{Lemma}

\theoremstyle{definition}
\newtheorem*{rem}{Remark}
\newtheorem{dfn}[theorem]{Definition}

\begin{document}

\author{Paolo Ghiggini} 
\address{CIRGET, Universit\'e du Qu\'ebec \`a Montr\'eal\\
Case Postale 8888, succursale Centre-Ville\\Montr\'eal (Qu\'ebec) H3C 3P8,
Canada}
\asciiaddress{CIRGET, Universite du Quebec a Montreal\\
Case Postale 8888, succursale Centre-Ville\\Montreal (Quebec) H3C 3P8,
Canada}

\email{ghiggini@math.uqam.ca}
\title{Strongly fillable contact $3$--manifolds\\without Stein fillings} 
\asciititle{Strongly fillable contact 3-manifolds without Stein fillings} 

\asciikeywords{Contact structure, symplectically fillable, Stein fillable, 
Ozsvath-Szabo invariant} 

\begin{abstract}
We use the Ozsv\'ath--Szab\'o contact invariant to produce
 examples of strongly symplectically fillable contact 
$3$--manifolds which are not Stein fillable.
\end{abstract}

\asciiabstract{%
We use the Ozsvath-Szabo contact invariant to produce
examples of strongly symplectically fillable contact 
3-manifolds which are not Stein fillable.}

\maketitle

\section{Introduction}
There is a strong relationship between contact topology and 
symplectic topology due to the fact that contact structures 
provide natural boundary conditions for symplectic structures 
on manifolds with boundary. Given a contact manifold $(Y, \xi)$ 
and a symplectic manifold $(W, \omega)$ with $\partial W=Y$, we say that 
$(W, \omega)$ {\em fills} $(Y, \xi)$ if some compatibility conditions 
are satisfied.  Depending on how restricting these conditions 
are, there are several different notions of fillability. The 
most widely studied in the literature are weak or strong 
symplectic fillability and Stein fillability. 

In the following we will always assume $Y$ is an oriented 
$3$--manifolds and $\xi$ is oriented and positive. This means 
that $\xi$ is the kernel of a
 globally defined smooth $1$--form $\alpha$ on $Y$ such that $\alpha \land d 
\alpha$ is a volume form inducing the fixed orientation of $Y$.

\begin{dfn}
A contact manifold $(Y, \xi)$ is {\em weakly symplectically fillable} if
$Y$ is the boundary of a symplectic manifold $(W, \omega)$ with $\omega|_{\xi}>0$.
\end{dfn}
Since $\omega$ orients $W$ and $\xi$ orients $Y$, we also require that the 
orientation of $Y$ as boundary of $W$ coincides with the orientation 
induced by $\xi$.

\begin{dfn}
A contact manifold $(Y, \xi)$ is {\em strongly symplectically fillable} if
$Y$ is the boundary of a symplectic manifold $(W, \omega)$ and $\xi$ is the
 kernel of a smooth $1$--form $\alpha$ on $Y$ such that $\omega|_{Y}= d \alpha$.
\end{dfn}

\begin{dfn}
A {\em Stein manifold} is a complex manifold $(X, J)$ with a proper 
function $\varphi\co X \to [0, + \infty)$ such that $dJ^*(d \varphi)$ is a K\"ahler form on 
$X$.
\end{dfn}

\begin{dfn}
A contact manifold $(Y, \xi)$ is {\em Stein fillable} (or
{\em holomorphically fillable}) if $Y$ is the boundary of a domain
$W= \varphi^{-1}([0,t])$ in a Stein manifold $(X, J)$ for some regular value
 $t$ of $\varphi$, and $\xi$ is the field of the complex hyperplanes of 
$J|_{\partial W}$.
\end{dfn}

\begin{figure}[ht!]\anchor{smooth.fig}
\centering
\psfrag{0}{\footnotesize $0$}
\psfrag{a}{\footnotesize $-n-1$}
\includegraphics[width=6cm]{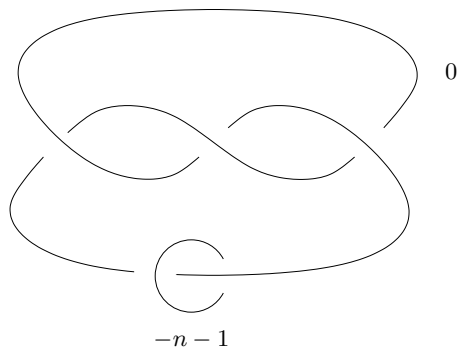}
\caption{The surgery diagram of $- \Sigma(2,3,6n+5)$}
\label{smooth.fig}
\end{figure}

\begin{rem}
In the literature there are several different equivalent definitions
of Stein manifold: see for example \cite[Section~4]{eliashberg-fill}.
\end{rem}
\def\strutt{\vrule width 0pt height 16pt depth 10pt}
There are obvious inclusions 
\[\left \{\strutt \parbox{1.35cm}{Stein \\ Fillable} \right \} \subset
\left \{\strutt \parbox{1.50cm}{Strongly \\ Fillable} \right \} \subset  
 \left \{\strutt \parbox{1.35cm}{Weakly \\ Fillable} \right \} \]
moreover, weakly fillable contact structures are tight by a deep 
theorem of Eliashberg and Gromov \cite{eliashberg:2,gromov}. The goal of
this article is to prove that the inclusion
\[\left \{\strutt \parbox{1.35cm}{Stein \\ Fillable} \right \} \subset
\left \{\strutt \parbox{1.50cm}{Strongly \\ Fillable} \right \}\]
is strict in dimension three. Let $- \Sigma(2,3,6n+5)$ be the $3$--manifold 
defined by the surgery diagram in \figref{smooth.fig}. We will 
prove the following theorem.

\begin{theorem}\label{main}
For any $n \geq 2$ and even, the $3$--manifold $- \Sigma(2,3,6n+5)$
 admits a strongly symplectically fillable contact structure 
which is not Stein fillable.
\end{theorem}

All other inclusions have already been proved to be strict: tight 
but non weakly fillable contact structures have been found first by 
Etnyre and Honda \cite{etnyre-honda:1} and later by Lisca and 
Stipsicz \cite{lisca-stipsicz:1,lisca-stipsicz:0}. A weakly fillable 
but non Strongly fillable contact structure has been found first by 
Eliashberg \cite{eliashberg:5} and later more have been found by 
Ding and Geiges \cite{ding-geiges:1}.

The main tool used in this article is the contact invariant in 
Heegaard--Floer theory recently introduced by Ozsv\'ath and Szab\'o 
\cite{O-Sz:cont}.  

\section{Construction of the non Stein fillable contact manifolds}
Let $M_0$ be the $3$--manifold obtained by $0$--surgery on the 
right-handed trefoil knot. $M_0$ admits  a presentation as 
a $T^2$--bundle over  $S^1$ with monodromy map
\[A\co T^2 \times \{ 1 \} \to T^2 \times \{ 0 \}\]
 given by 
$A= \left ( \begin{matrix}
                          1 & 1 \\
                         -1 & 0 
                                 \end{matrix} \right )$.
Put coordinates $(x,y,t)$ on $T^2 \times \R$. The $1$--forms 
\[\alpha_n = \sin(\phi(t))dx+ \cos(\phi(t))dy\]
on $T^2 \times \R$ define contact structures $\xi_n$ on $M_0$ for any $n>0$ 
provided that 
\begin{enumerate}
\item $\phi'(t)>0$ for any $t \in \R$
\item $\alpha_n$ is invariant under the action 
$({\mathbf v},t) \mapsto (A {\mathbf v}, t-1)$
\item $n \pi \leq \sup \limits_{t \in \R} (\phi(t+1)- \phi(t))< (n+1) \pi$.
\end{enumerate}

The main result about this family of contact structures we will need
 in the present article is the following.
\begin{theorem}\label{giga}{\rm(\cite[Proposition~2 and Theorem~6]{giroux:3}, 
\cite[Theorem~1]{ding-geiges:1})}\qua
The contact structures $\xi_n$ do not depend on the function $\phi$ up to 
isotopy, and are all weakly symplectically fillable.
\end{theorem}

 Let $F$ be the image in $M_0$ of the segment $\{ 0 \} \times [0,1] \subset 
T^2 \times [0,1]$, then $F$ is Legendrian 
with respect to the contact structure $\xi_n$ for all $n$. Denote by $K$
the right-handed trefoil knot in $S^3$. We can 
choose a diffeomorphism from the complement of a tubular 
neighbourhood of  $K$ in $S^3$ to the 
complement of a tubular neighbourhood of $F$ in $M_0$ so that the 
meridian of $K$ is mapped to a longitude of $F$. This diffeomorphism 
defines a framing on $F$, and the 
framing so defined allows us to define a twisting number for $F$.

\begin{lemma}\label{twisting}{\rm\cite[Lemma~3.5]{ghiggini:3}}\qua
The twisting number of $\xi_n$ along the Legendrian curve $F$ is
$tn(F, \xi_n)= -n$
\end{lemma}
Legendrian surgery on $(M_0, \xi_n)$ along $F$ is smoothly equivalent to
the surgery described by \figref{smooth.fig} which produces the 
manifold $-\Sigma(2,3,6n+5)$.
We denote the tight contact structure on $-\Sigma(2,3,6n+5)$ obtained by 
Legendrian surgery on $(M_0, \xi_n)$ along $F$ by $\eta_0$. The following 
theorem proves the strong fillability part of Theorem \ref{main}.

\begin{theorem}
The contact manifolds $(-\Sigma(2,3,6n+5), \eta_0)$ are strongly 
symplectically fillable for any $n \geq 1$.
\end{theorem}
\begin{proof}
The contact manifolds $(M_0, \xi_n)$ are weakly symplectically fillable 
by Theorem \ref{giga}.  Since Legendrian surgery 
preserves weak fillability by \cite[Theorem~2.3]{etnyre-honda:2},
$(- \Sigma(2,3,6n+5), \eta_0)$ is also weakly fillable.
Since the manifolds $\Sigma(2,3,6n+5)$ are
 homology spheres, by  \cite[Proposition~4.1]{eliashberg-fill} the 
symplectic form on the filling can be modified in a neighbourhood 
of the boundary so that the filling becomes strong.
\end{proof} 
The non Stein fillability part of Theorem \ref{main} can now be made 
more precise with the following 
statement.
\begin{theorem}\label{nonstein}
The contact manifolds $(-\Sigma(2,3,6n+5), \eta_0)$ are not Stein fillable 
for any $n \geq 2$ and even.
\end{theorem}
The proof of this theorem is the goal of Section 4.

\section{Overview of the contact invariant}
In this section we give a brief overview of Heegaard--Floer homology 
and of the contact invariant defined by Ozsv\'ath and Szab\'o. We will 
not treat the subject in its most general form, but only in the 
form it will be needed in the proof of Theorem \ref{nonstein}. 

\subsection{Heegaard--Floer homology}
Heegaard--Floer homology is a family of topological quantum field 
theories for $Spin^c$ $3$--manifolds introduced by Ozsv\'ath and Szab\'o 
in \cite{O-Sz:3,O-Sz:2,O-Sz:1}. 
In their simpler form they associate vector spaces $\widehat{HF}(Y, 
\mathfrak{t})$ and $HF^+(Y, \mathfrak{t})$ over $\Z /2 \Z$ to any 
closed oriented $Spin^c$ $3$--manifold $(Y, \mathfrak{t})$, and  
homomorphisms 
\[F^{\circ}_{W, \mathfrak{s}}\co HF^{\circ}(Y_1, \mathfrak{t_1}) \to 
HF^{\circ}(Y_2, \mathfrak{t_2})\]
to any oriented $Spin^c$--cobordism $(W, \mathfrak{s})$ between two 
$Spin^c$--manifolds $(Y_1, \mathfrak{t_1})$ and $(Y_2, \mathfrak{t_2})$ 
such that $\mathfrak{s}|_{Y_i}= \mathfrak{t_i}$. Here $HF^{\circ}$ denotes 
either $\widehat{HF}$ or $HF^+$. The groups $\widehat{HF}(Y, 
\mathfrak{t})$ and $HF^+(Y, \mathfrak{t})$ are linked to one another 
by the exact triangle
\begin{align} \label{triangolo}
\parbox{12cm}{\xymatrix{ 
& \widehat{HF}(Y, \mathfrak{t}) \ar[r] & HF^+(Y, \mathfrak{t}) \ar[r] 
& HF^+(Y, \mathfrak{t}) \ar `r[d] `[l] `[lll] `[l] [ll] \\
& & &       
}}
\end{align}
This exact triangle is natural in the sense that its maps commute 
with the maps induced by cobordisms.

It was shown in \cite{O-Sz:3} that, when $c_1(\mathfrak{t})$ is a 
torsion element, the vector spaces $\widehat{HF}(Y, \mathfrak{t})$ 
and $HF^+(Y, \mathfrak{t})$ come with a $\Q$--grading.  In conclusion, 
for a torsion $Spin^c$--structure $\mathfrak{t}$ on $Y$ the 
Heegaard--Floer homology groups $HF^{\circ} (Y, \mathfrak{t})$ split as
\[HF^{\circ} (Y, \mathfrak{t})=\bigoplus_{d \in  \Q} HF^{\circ}_{(d)} (Y, \mathfrak{t}).\] 

The set of the $Spin^c$--structures on a manifold has an involution 
called {\em conjugation}. Given a $Spin^c$--structure $\mathfrak{t}$, we 
denote its conjugate $Spin^c$--structure by $\overline{\mathfrak{t}}$. 
We have $c_1(\overline{\mathfrak{t}})= -c_1(\mathfrak{t})$. 
There is an isomorphism $\mathfrak{J}\co HF^{\circ}(Y, \mathfrak{t}) \to 
HF^{\circ}(Y, \overline{\mathfrak{t}})$ defined in \cite[Theorem~2.4]{O-Sz:2}. 
We recall that the isomorphism $\mathfrak{J}$ preserves 
the $\Q$--grading of the Heegaard--Floer homology groups when 
$c_1(\mathfrak{t})$ is a torsion cohomology class, and is a natural 
transformation in the following sense.
\begin{prop}{\rm\cite[Theorem~3.6]{O-Sz:3}}\qua
Let $(W, \mathfrak{s})$ be a $Spin^c$--cobordism between 
$(Y_1, \mathfrak{t}_1)$ and $(Y_2, \mathfrak{t}_2)$. Then  
the following  diagram 
$$\begin{CD}
HF^{\circ}(Y_1, \mathfrak{t}_1) & @> F^{\circ}_{W, \mathfrak{s}} >> & HF^{\circ}(Y_2, \mathfrak{t}_2) \\
@VV \mathfrak{J} V &   & @VV \mathfrak{J} V \\
HF^{\circ}(Y_1, \overline{\mathfrak{t}}_1) & @> F^{\circ}_{W, \overline{\mathfrak{s}}} >> &
 HF^{\circ}(Y_2, \overline{\mathfrak{t}}_2)
\end{CD}$$
commutes. 
\end{prop}
The isomorphism $\mathfrak{J}$ commutes also with the maps in the 
exact triangle (\ref{triangolo}).

\subsection{Contact invariant}
A contact structure $\xi$ on a $3$--manifold $Y$ determines a 
$Spin^c$--structure
$\mathfrak{t}_{\xi}$ on $Y$ such that $c_1(\mathfrak{t}_{\xi})=c_1(\xi)$.
To any contact manifold $(Y, \xi)$ we can associate an element $c(\xi) \in 
\widehat{HF}(-Y, \mathfrak{t}_{\xi})$ which is an isotopy invariant of $\xi$,
see \cite{O-Sz:cont}.
Sometimes it is also useful to consider the image $c^+(\xi) \in HF^+(-Y, 
\mathfrak{t}_{\xi})$ of $c(\xi)$ under the  map $\widehat{HF}(-Y, 
\mathfrak{t}_{\xi}) \to HF^+(-Y, \mathfrak{t}_{\xi})$ in the exact triangle 
(\ref{triangolo}). The Ozsv\'ath--Szab\'o contact 
invariant satisfies the following properties.
\begin{theorem} {\rm\cite[Theorem~1.4 and Theorem~1.5]{O-Sz:cont}}\qua
If $(Y, \xi)$ is overtwisted, then $c(\xi)=0$. If $(Y, \xi)$ is Stein fillable, 
then $c(\xi) \neq 0$.
\end{theorem}
\begin{prop}{\rm\cite[Proposition~4.6]{O-Sz:cont}}\qua
If $c_1(\xi)$ is a torsion homology class, then $c(\xi)$ is a homogeneous 
element of degree $-d_3(\xi)- \frac 12$, where $d_3(\xi)$ denotes the 
$3$--dimensional homotopy invariant introduced by Gompf 
\cite[Definition~4.2]{gompf:1}. 
\end{prop}
\begin{theorem}\label{olga1}{\rm\cite[Theorem~4]{plam:1}}\qua
Let $W$ be a smooth compact $4$--manifold with boundary $Y= \partial W$. Let
$J_1$, $J_2$ be two Stein structures on $W$ that induce $Spin^c$--structures 
$\mathfrak{s}_1$, $\mathfrak{s}_2$ on $W$ and
contact structures $\xi_1$, $\xi_2$ on $Y$. We puncture $W$ and regard it as a 
cobordism from $-Y$ to $S^3$. Suppose that $\mathfrak{s}_1|_Y$ is isotopic to
$\mathfrak{s}_2|_Y$, but the $Spin^c$--structures $\mathfrak{s}_1$, 
$\mathfrak{s}_2$ are not isomorphic. Then
\begin{enumerate}
\item $F^+_{W, \mathfrak{s}_i}(c(\xi_j)) = 0$ for $i \neq j$;
\item $F^+_{W, \mathfrak{s}_i}(c(\xi_i))$ is a generator of $HF^+(S^3)$.
\end{enumerate}
\end{theorem}

The space of oriented contact structures on $Y$ has a natural involution
 called {\em conjugation}. For any contact structure $\xi$ on a 
$3$--manifold $Y$ we denote by $\overline{\xi}$ the contact structure on 
$Y$ obtained from $\xi$ by inverting the orientation of the planes. 
The conjugation of contact structures is compatible with  the 
conjugation of the $Spin^c$--structure defined by the contact structure, 
in fact $\mathfrak{t}_{\overline{\xi}} = \overline{\mathfrak{t}_{\xi}}$. The contact
 invariant behaves well with respect to conjugation.
\begin{prop}\label{mio}{\rm\cite[Theorem~2.10]{ghiggini:3}}\qua
Let $(Y, \xi)$ be a contact manifold, then 
\[c(\overline{\xi})= \mathfrak{J}(c(\xi)).\]
\end{prop}

\section{Proof of the  non fillability of $(- \Sigma(2,3,6n+5), \eta_0)$}
In this article we will consider only integer homology spheres,  
which have therefore a unique $Spin^c$--structure. For this reason from 
now on we will always suppress the $Spin^c$--structure in the notation of 
the Heegaard--Floer groups.

The key ingredients in the proof of Theorem \ref{nonstein} are the 
conjugation invariance of $\eta_0$ and the structure of the 
$\mathfrak{J}$--action on $\widehat{HF}(\Sigma(2,3,6n+5))$. The starting 
point is a general 
observation about the Stein fillings of conjugation invariant 
contact structures.
\begin{prop}\label{spin}
Let $\xi$ be a contact structure on a $3$--manifold $Y$ which is isotopic 
to its conjugate $\overline{\xi}$. If $(W, J)$ is a Stein filling of $\xi$ 
and $\mathfrak{s}$ is its canonical $Spin^c$--structure, then 
$\mathfrak{s}$ is isomorphic to its conjugate $\overline{\mathfrak{s}}$.
\end{prop}
\begin{proof}
If $(W, J)$ is a Stein filling of $\xi$, then $(W, -J)$ is a Stein 
filling of $\overline{\xi}$, and the canonical $Spin^c$--structure of 
$(W, -J)$ is $\overline{\mathfrak{s}}$. Puncture $W$ and regard it as a 
cobordism between $-Y$ and $S^3$. Since $\overline{\xi}$ is isotopic to 
$\xi$ we have 
\[F_{W, \mathfrak{s}}(c(\xi))= F_{W, \overline{\mathfrak{s}}}(c(\xi)) \neq 0.\] 
Theorem \ref{olga1} implies that $\mathfrak{s}$ is isomorphic to 
$\overline{\mathfrak{s}}$. 
\end{proof}

\begin{rem}
Proposition \ref{spin} can be deduced also from Seiberg--Witen theory, see for 
example \cite[Theorem~1.2]{lisca-matic:1} or 
\cite[Theorem~1.2]{kronheimer-mrowka:1}.
\end{rem}

By \cite[Theorem~3.12]{ghiggini:3}  the $3$--dimensional homotopy
invariant of $\eta_0$ is $d_3(\eta_0)=- \frac 32$, therefore the contact invariant
$c(\eta_0)$ belongs to $\widehat{HF}_{(+1)}(\Sigma(2,3,6n+5))$. The group
$HF^+(- \Sigma(2,3,6n+5))$ is computed in \cite[Section~8]{O-Sz:4}. From this
it is easy to prove that $\widehat{HF}_{(+1)}(\Sigma(2,3,6n+5))$ is 
isomorphic to $(\Z / 2 \Z)^n$ by applying the exact triangle 
(\ref{triangolo})
and the isomorphism $\widehat{HF}_{(d)}(Y) \cong \widehat{HF}_{(-d)}(-Y)$
which holds for any homology sphere $Y$.

\begin{figure}[ht!]\anchor{legendrian.fig}
\centering
\psfrag{a}{\footnotesize $\frac{n-i}{2}$ cusps}
\psfrag{b}{\footnotesize $\frac{n+i}{2}$ cusps}
\includegraphics[width=6cm]{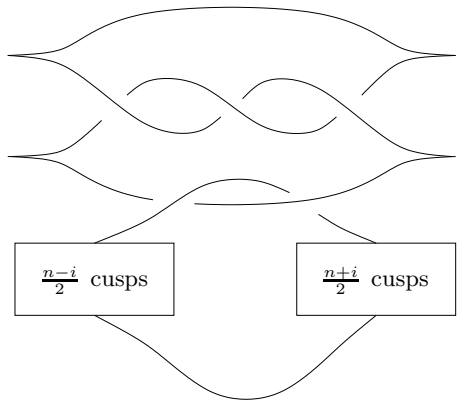}
\caption{Legendrian surgery presentation of the contact manifold
$(- \Sigma(2,3,6n+5), \eta_i)$ for $i \in {\mathcal P}^*_n$}
\label{legendrian.fig}
\end{figure}

Now we give a closer look at the action of $\mathfrak{J}$ on 
$\widehat{HF}_{(+1)}(\Sigma(2,3,6n+5))$ by considering the action of conjugation 
on a set of Stein fillable contact structures on $- \Sigma(2,3,6n+5)$.
For any $n \in \N$ and $n \geq 2$ we define 
\[{\mathcal P}^*_n = \{ -n+1, -n+3, \ldots ,n-3, n-1 \}.\] 
If $n$ is even, then $0 \notin {\mathcal P}^*_n$.  
Given $i \in {\mathcal P}^*_n$, by $\eta_i$ we denote the contact 
structure on $-\Sigma(2,3,6n+5)$ obtained by Legendrian surgery on 
the Legendrian link in the standard $S^3$ shown in
\figref{legendrian.fig}. In the following we will 
always assume $n$ even, so there is no confusion between $\eta_0$ as 
defined in Section 2 and $\eta_i$ with $i \in {\mathcal P}^*_n$.
 The contact structures $\eta_i$ with $i \in {\mathcal P}^*_n$ are all 
Stein fillable and pairwise homotopic with $3$--dimensional 
homotopy invariant $d_3(\eta_i)=- \frac 32$. 

\begin{prop} \label{olga2} {\rm\cite[Section~4]{plam:1}}\qua
The contact invariants $c(\eta_i)$ for $i \in {\mathcal P}^*_n$ 
generate $\widehat{HF}_{(+1)}(\Sigma(2,3,6n+5))$.
\end{prop}

\begin{prop} \label{coniugio}{\rm\cite[Proposition~3.8]{ghiggini:3}}\qua
The contact structure $\overline{\eta_i}$ obtained from $\eta_i$ by 
conjugation is isotopic to $\eta_{-i}$ for $i \in {\mathcal P}^*_n$, 
and $\eta_0$ is isotopic to its conjugate $\overline{\eta_0}$.
\end{prop}
Putting Proposition \ref{olga2} and Proposition \ref{coniugio} 
together we obtain the following lemma.
\begin{lemma} 
If $n$ is even, then the subspace 
\[\mbox{Fix}(\mathfrak{J}) \subset \widehat{HF}_{(+1)}(\Sigma(2,3,6n+5))\] 
of the fix points for the action of $\mathfrak{J}$ on 
$\widehat{HF}_{(+1)}(\Sigma(2,3,6n+5))$ is generated by 
$c(\eta_i)+c(\eta_{-i})$ for $i \in {\mathcal P}^*_n$. 
\end{lemma}
\begin{proof}
Let $x \in \mbox{Fix}(\mathfrak{J})$ be a fixed point. We write 
 \[ x= \sum_{i \in {\mathcal P}^*_n} \alpha_i c(\eta_i) \]
for $\alpha_i \in \{ 0,1 \}$, then applying $\mathfrak{J}$ we obtain
\[ x= \sum_{i \in {\mathcal P}^*_n} \alpha_i c(\eta_{-i}). \]
From this we deduce that $\alpha_i = \alpha_{-i}$, which implies the lemma.
\end{proof}

\begin{proof}[Proof of Theorem \ref{nonstein}]
 Suppose  $(W, J)$ is a 
Stein filling of $(- \Sigma(2,3,6n+5), \eta_0)$ and call 
$\mathfrak{s}$ its canonical $Spin^c$-structure. By 
Proposition \ref{spin} $\mathfrak{s}$ is invariant under 
conjugation. Moreover, $c(\eta_0) \in \mbox{Fix}(\mathfrak{J})$ by 
Proposition \ref{mio}, therefore $c(\eta_0)$ is a linear 
combination of elements of the form $c(\eta_i)+c(\eta_{-i})$ for 
$i \in {\mathcal P}^*_n$.
Applying the map $\widehat{HF}(\Sigma(2,3,6n+5)) \to  
HF^+(\Sigma(2,3,6n+5))$ we obtain that $c^+(\eta_0)$ is a linear 
combination of elements of the form $c^+(\eta_i)+c^+(\eta_{-i})$. 

Puncture the Stein filling $W$ and regard it as a cobordism from 
$\Sigma(2,3,6n+5)$ to $S^3$. Applying $F^+_{W, \mathfrak{s}}$ to each 
$c^+(\eta_i)+c^+(\eta_{-i})$ we get
\[F^+_{W, \mathfrak{s}}(c^+(\eta_i)+c^+(\eta_{-i}))= 
F^+_{W, \mathfrak{s}}(c^+(\eta_i))+ F^+_{W,
  \mathfrak{s}}(\mathfrak{J}(c^+(\eta_i)) 
= 2F^+_{W, \mathfrak{s}}(c^+(\eta_i))=0\]
because 
\[F^+_{W, \mathfrak{s}}(\mathfrak{J}(c^+(\eta_i))= 
\mathfrak{J}(F^+_{W, \overline{\mathfrak{s}}}(c^+(\eta_i)))=F^+_{W, \mathfrak{s}}(c^+(\eta_i))\]
by Proposition \ref{spin}, the naturality of the homomorphism 
$\mathfrak{J}$, and the triviality of the $\mathfrak{J}$--action
on $HF^+(S^3)$.
This implies $F^+_{W, \mathfrak{s}}(\eta_0)=0$, which is a contradiction with 
Theorem \ref{olga1}(2), therefore a Stein filling of 
$(- \Sigma(2,3,6n+5), \eta_0)$ 
cannot exist. 
\end{proof}

\begin{rem}
With the same argument we can actually prove that $(- \Sigma(2,3,6n+5), \eta_0)$
has no symplectic filling with exact symplectic form when $n$ is even.
We will call such a filling an {\em exact filling}. Exact fillability
is a notion of fillability which is intermediate between strong and 
Stein fillability and has not been studied much yet. We do not know at 
present if exact fillability is a 
different notion from Stein fillability.

This stronger form of Theorem \ref{nonstein} can be proved  by extending
Theorem \ref{olga1} to exact fillings.
\end{rem}

\end{document}